\magnification=\magstephalf
\input eplain
\input BMmacs

\phantom{x}\vskip .5cm
\centerline{{\titlefont ON A QUESTION OF O'GRADY ABOUT MODIFIED DIAGONALS}}
\vskip .5cm

\centerline{{\byfont by}\quad {\namefont Ben Moonen and Qizheng Yin}}
\vskip 1cm

{\eightpoint 
\noindent
{\bf Abstract.} Let $X$ be an abelian variety of dimension~$g$. In a recent preprint O'Grady defines modified diagonal classes $\Gamma^m$ on~$X^m$ and he conjectures that the class of~$\Gamma^m$ in the Chow ring of~$X^m$ is torsion for $m \geq 2g+1$. We prove a generalization of this conjecture.
\medskip

\par}
\vskip 8mm

\section
\sectlabel{Notat}
Throughout this note, $R$ is a Dedekind ring and $S$ is a connected scheme that is smooth and of finite type over~$R$, which is the base scheme over which we work.

Let $\pi \colon Y \to S$ be a smooth projective $S$-scheme, $a \colon S \to Y$ a section of~$\pi$, and $m \geq 1$ an integer. By $Y^m$ we mean the product of $m$~copies of~$Y$ relative to~$S$. For $j \subset \{1,\ldots,m\}$ let $\pr_j \colon Y^m \to Y$ be the projection onto the $j$th factor and let $\hat\pr_j \colon Y^m \to Y^{m-1}$ be the projection that contracts the $j$th factor.

For $I \subset \{1,\ldots,m\}$, let $\delta_I \colon Y \to Y^m$ be the morphism defined by the property that $\pr_i \circ \delta_I$ is the identity on~$Y$ if $i \in I$ and is the ``constant map'' $a \circ \pi$ if $i \notin I$. Following O'Grady~[\ref{OGr}] we define an algebraic cycle $\Gamma^m(Y,a)$ on~$Y^m$ by
$$
\Gamma^m(Y,a) = \sum_{\emptyset \neq I \subset \{1,\ldots,m\}} (-1)^{m-|I|} \cdot \delta_{I,*}(Y)\, .
$$

The goal of this note is to prove the following result, which proves Conjecture~0.3 of~[\ref{OGr}], generalized to the relative setting.

\section
\sectlabel{MainThm}
{\it Theorem. --- Let $X \to S$ be an abelian scheme of relative dimension~$g$. For any section $a \in X(S)$ the class of~$\Gamma^m(X,a)$ in $\CH(X^m)$ is torsion if $m \geq 2g+1$.
\par}

\section
{\it Remark.\/} Still with $X$ an abelian scheme over~$S$, let $e \in X(S)$ be the zero section. Translation over~$a$ gives an isomorphism $t_a \colon X \isomarrow X$ over~$S$. The class of $\Gamma^m(X,a)$ is the push-forward of the class of~$\Gamma^m(X,e)$ under $t_a^m \colon X^m \isomarrow X^m$. Hence it suffices to prove the theorem taking~$e$ as section.

\section
\sectlabel{PushFLem}
{\it Lemma. --- Let $X/S$ be an abelian scheme of relative dimension~$g$ and write $\Gamma^m = \Gamma^m(X,e)$. For $n \in \mZ$ let $\mult_X(n) \colon X \to X$ be the endomorphism given by multiplication by~$n$.

\item{\rm (\romno1)} For $n \in \mZ$ we have $\mult_{X^m}(n)_* [\Gamma^m] = n^{2g} \cdot [\Gamma^m]$.

\item{\rm (\romno2)} For $j \in \{1,\ldots,m\}$ let $\hat\pr_j \colon X^m \to X^{m-1}$ be the projection map contracting the $j$th factor. Then $\hat\pr_{j,*}[\Gamma^m] = 0$.
\par}
\medskip

\Proof
For (\romno1) we use that, with notation as in point~\ref{Notat}, and taking $a=e$ as section, $\mult_{X^m}(n) \circ \delta_I$ is the same as $\delta_I \circ \mult_X(n)$ and that $\mult_X(n)$ is finite flat of degree~$n^{2g}$.

For (\romno2), consider a non-empty subset $I \subset \{1,\ldots,m\}$. If $I = \{j\}$ then $\hat\pr_j \circ \delta_I \colon X \to X^{m-1}$ is constant, so that $\hat\pr_{j,*} \delta_{I,*}(X) = 0$.  Let $\cI$ be the set of non-empty subsets of $\{1,\ldots,m\}$ different from~$\{j\}$. Then $\cI = \cI_0 \coprod \cI_1$ where $I \in \cI_0$ if $j \notin I$ and $I \in \cI_1$ if $j \in I$. The map $\beta\colon \cI_0 \to \cI_1$ given by $I \mapsto I\cup \{j\}$ is a bijection. Further, for $I \in \cI_0$ we have $\hat\pr_j \circ \delta_I = \hat\pr_j \circ \delta_{\beta(I)}$ and $|\beta(I)| = |I|+1$. In the calculation of $\hat\pr_{j,*} [\Gamma^m]$, the terms corresponding to~$I$ and~$\beta(I)$ therefore cancel, and this gives the assertion.
\QED

\section
Let $\Mot_S$ be the category of Chow motives over~$S$ with respect to graded correspondences; see for instance [\ref{DenMur}], Section~1. If $Y$ is a smooth projective $S$-scheme we write $h(Y)$ for its Chow motive. In $\Mot_S$ we have a tensor product with $h(Y) \otimes h(Z) = h(Y\times_S Z)$.

Let $\unitmot_S = h(S)$ be the unit motive and $\unitmot(n)$ the $n$th Tate twist. If $M$ is a motive, we write $M(n) = M \otimes \unitmot(n)$. The Chow groups (with $\mQ$-coefficients) of a motive~$M$ are defined by $\CH^i(M)_\mQ = \Hom_{\Mot_S}\bigl(\unitmot(-i),M\bigr)$. 

If $f\colon Y \to Z$ is a morphism of smooth projective $S$-schemes we have induced morphisms $f^* \colon h(Z) \to h(Y)$ and, assuming $Y$ and~$Z$ are connected, $f_*\colon h(Y) \to h(Z)\bigl(d\bigr)$, where $d = \dim(Z/S) - \dim(Y/S)$.

\section
Let $X/S$ be an abelian scheme of relative dimension~$g$. As proven by Deninger and Murre in~[\ref{DenMur}] (generalizing results of Beauville~[\ref{BeauvChow}] over a field) we have a canonical decomposition $h(X) = \oplus_{i=0}^{2g}\,  h^i(X)$ in~$\Mot_S$ that is stable under all endomorphisms~$\mult_X(n)_*$, and such that $\mult_X(n)_*$ is multiplication by~$n^{2g-i}$ on~$h^i(X)$. For $m \geq 1$ this induces a decomposition 
$$
h(X^m) = \bigoplus_{\bi = (i_1,\ldots,i_m)}\quad \bigotimes_{j=1}^m\; h^{i_j}(X) \, ,  
$$
where the sum runs over the elements $\bi \in \{0,\ldots,2g\}^m$. Under this decomposition we have
$$
h^\nu(X^m) = \bigoplus_{|\bi|=\nu}\quad \bigotimes_{j=1}^m\; h^{i_j}(X) \, ,  \eqlabel{hnuXmDec}
$$
where the sum runs over the $m$-tuples $\bi = (i_1,\ldots,i_m)$ in $\{0,\ldots,2g\}^m$ with $|\bi| = i_1+\cdots+i_m$ equal to~$\nu$.

If $\pi \colon X \to S$ is the structural morphism, $\pi_* \colon h(X) \to h(S)\bigl(-g\bigr) = \unitmot(-g)$ is an isomorphism on~$h^{2g}(X)$ and is zero on $\oplus_{i=0}^{2g-1}\,  h^i(X)$

\section
\sectlabel{VanishLem}
{\it Lemma. ---  Notation as above. If there is an index~$\nu$ such that $i_\nu = 2g$ then the component of $\bigl[\Gamma^m(X,e)\bigr]$ in $\CH(\otimes_{j=1}^m\, h^{i_j}(X))_\mQ$ is zero.
\par}
\medskip

\Proof
Assume $i_\nu = 2g$. Consider the projection $\hat\pr_\nu \colon X^m \to X^{m-1}$. By the fact stated just before the lemma, the induced map $\hat\pr_{\nu,*} \colon h(X^m) \to h(X^{m-1})\bigl(-g\bigr)$ restricts to an isomorphism of $\otimes_{j=1}^m\, h^{i_j}(X)$ with a sub-motive of $h(X^{m-1})\bigl(-g\bigr)$. The assertion therefore follows from Lemma~\ref{PushFLem}(\romno2).
\QED

\section
{\it Proof of the Theorem.\/} Write $\Gamma^m = \Gamma^m(X,e)$. By Lemma~\ref{PushFLem}(\romno1) we have
$$
[\Gamma^m] \in \CH\bigl(h^{2g(m-1)}(X^m)\bigr)_\mQ \subset \CH(X^m)_\mQ\, ,
$$
and \eqref{hnuXmDec} gives
$$
\CH\bigl(h^{2g(m-1)}(X^m)\bigr)_\mQ = \bigoplus_{|\bi|=2g(m-1)}\quad \CH(\otimes_{j=1}^m\, h^{i_j}(X))_\mQ\, .
$$
The assumption that $m \geq 2g+1$ implies that for every $\bi = (i_1,\ldots,i_m) \in \{0,\ldots,2g\}^m$ with $|\bi| = 2g(m-1)$ there is an index~$\nu$ with $i_\nu = 2g$ and by Lemma~\ref{VanishLem} we are done. \QED

\vskip2.0\bigskipamount plus 2pt minus 1pt%
\goodbreak\noindent{{\bf References}}%
\nobreak\vskip.75\bigskipamount plus 2pt minus 1pt%

{\eightpoint

\bibitem{BeauvChow}
A.~Beauville, {\it Sur l'anneau de Chow d'une vari\'et\'e ab\'elienne.\/}
Math.\ Ann.\ 273 (1986), 647--651. 

\bibitem{DenMur}
C.~Deninger, J.~Murre, {\it Motivic decomposition of abelian schemes and the Fourier transform.\/} J.\ reine angew.\ Math.\ 422 (1991), 201--219.

\bibitem{OGr}
K.G.\ O'Grady, {\it Computations with modified diagonals}. Preprint, arXiv:1311.0757v1.

\vskip 1cm

\noindent
Radboud University Nijmegen, IMAPP, PO Box 9010, 6500GL Nijmegen, The Netherlands. 
\smallskip

\noindent
b.moonen@science.ru.nl\qquad q.yin@science.ru.nl
\par}

\bye